\documentclass[12pt]{article}
\usepackage[english]{babel}
\usepackage[a4paper, total={6.5in, 9in}]{geometry}
\usepackage{graphicx} 
\usepackage{amsmath}
\usepackage{amsfonts}
\usepackage{faktor}
\usepackage{amssymb}
\usepackage{placeins}
\usepackage[hidelinks]{hyperref}
\usepackage{comment}
\usepackage{tikz}
\usepackage{csquotes}
\usepackage[backend=biber]{biblatex}
\usepackage[T1]{fontenc}

\newcommand{\dron}[1]{\partial_{#1}}
\newcommand{\mc}[1]{\mathcal{#1}}
\newcommand{\ep}{\varepsilon}
\newcommand{\R}{\mathbb{R}}

\newcommand{\ZZ}{\mathbb{Z}}
\newcommand{\fee}{\varphi}
\newcommand{\map}{\xrightarrow{}}
\newcommand{\inv}{^{-1}}

\newcommand{\SSS}{\mathcal{S}}
\newcommand{\Lp}{\widehat L^+}
\newcommand{\Lm}{\widehat L^-}
\newcommand{\Kp}{\widehat K^+}

\newcommand{\Km}{\widehat K^-}
\newcommand{\Kch}{\widehat K}
\newcommand{\Lch}{\widehat L}
\newcommand{\SingZ}{\{0\} }
\newcommand{\reeb}{\textbf{Reeb} }
\newcommand{\nrmt}{\|t\| }
\newcommand{\nrmtc}{\|t\|^2}
\newcommand{\lbar}{{\overline{L}}}
\newcommand{\feetild}{\tilde{\varphi}}

\usepackage{amsthm}
\newtheorem{thm}{Theorem}

\newtheorem{deff}[thm]{Definition}
\newtheorem{prop}[thm]{Proposition}
\newtheorem{lemm}[thm]{Lemma}
\newtheorem{rmq}[thm]{Remark}
\newtheorem*{findepreuve}{Proof of theorem 1}

\title{Lagrangian Intersections in Liouville Manifolds via Generating Functions}
\author{Yvon Fredrich}

\begin{document}
\maketitle

\section{Context and introduction}

Let $(W, \lambda)$ be a Liouville manifold and consider a Legendrian isotopy $(L_t)_{t\in[0, 1]}\subset \mc C W = W\times\R$  where $(\mc C W, dz-\lambda)$ is the contactisation of $W$. We seek a lower bound for $\#\reeb(L_0, L_1)$ where $$\textbf{Reeb}(L_0, L_1):=\{ (x, z_0, z_1)\in W\times\R\times\R \,|\, (x, z_0)\in L_0, (x, z_1)\in L_1\}.$$
In the particular case where $W$ is a cotangent bundle, analogous results are proved via the theory of generating families by Laudenbach and Sikorav in \cite{Laudenbach_1985} and by Chekanov in \cite{Chekanov_1996}. Laudenbach also proved such a result for a general Liouville manifold $W$ in \cite{Laudenbach_1995}, with the additional hypothesis $\pi_1(W, L) = 0, \pi_2(W, L) = 0$.

\subsection{General outline of the proof}
We embed $W$ symplectically in a cotangent bundle  $T^*M$, and construct a legendrian  $\Lch \subset J^1M$ (The $1$-jet space on $M$) via a generating family $g$ (the formula for $g$ was introduced in \cite{abouzaid2024twistedgeneratingfunctionsnearby}). The Legendrian $\Lch$ is a sphere bundle over $L$, with a Reeb chord between the north and south poles of the sphere above each point of $L$. The idea of embedding W into a cotangent bundle  was introduced by Shende in the context of microlocal sheaf theory in  \cite{Shende2017} (see also \cite{nadler2022sheafquantizationweinsteinsymplectic} and \cite{Eliashberg1998-mz} section $0.3.9$). We use the isotopy $L_s$ to obtain an isotopy $\Lch_s$ such that the long Reeb chords of  $\Lch_1$ are in bijection with those between $L_0$ and $L_1.$ Chekanov's theorem \cite{Chekanov_1996} provides an isotopy of generating families $g_s$ with $L_{g_s} = \Lch_s$ and the result then follows from studying the sublevel homology of $\Lch_s.$

\begin{thm}
    Let $(W, \lambda)$ be a Liouville manifold equipped with a stable polarisation $P$, and $(L_t)_{t\in[0, 1]}\subset \mc CW$ a legendrian isotopy. Suppose $L_0$ is compact and $\pi_{|L_0}$ is an embedding, where $\pi : \mc C W\map W$ is the projection onto the first factor. Suppose also that $T\pi(L_0)$  is stably homotopic to $P_{|\pi(L_0)}$ as stable lagrangian subbundles of  $TW_{|\pi(L_0)}.$
    If  $\pi(L_0) \pitchfork \pi(L_1) $ then $$\#\pi(L_0)\cap \pi(L_1) = \# \reeb(L_0, L_1)\geq \displaystyle\sum_i b_i(L_0).$$
    If also $L_0\hookrightarrow W$ induces a left invertible map $\pi_1(L_0)\map \pi_1(W)$ then
    $$\#\pi(L_0)\cap \pi(L_1) = \# \reeb(L_0, L_1)\geq \mu_{st}(L),$$
    where $\mu_{st}$ is the stable Morse number, defined below.
    
\end{thm}

\begin{rmq}
    By slightly tweaking the proof of lemma $(12)$ we can weaken the hypothesis and suppose only that $L_0$ be  hamiltonian isotopic to a legendrian $L$ satisfying that $\pi_{|L}$ is an embedding and $T\pi(L)$  stably homotopic to $P_{|\pi(L)}$ as stable lagrangian subbundles of  $TW_{|\pi(L)}.$ 
\end{rmq}

The purpose of the following proof is mainly to generalise the use of generating families to Liouville manifolds. Indeed, supposing some results of legendrian contact homology we can give a very short proof sketch :
Let $\phi_R^s : \mc CW\map \mc CW$ be the Reeb flow. For $T$ sufficiently large, $\displaystyle\min_{t\in [0, 1]}\displaystyle\min_{x \in L_t} z(\phi_R^T(x)) > \displaystyle\max_{x\in L_0} z(L_0)$. 
We obtain an isotopy of legendrian links : $(L_0, \phi_R(L_t))_t$. Starting from the trivial augmentation  $\ep_0'$ of $\phi_R(L_0)$ we obtain an augmentation $\ep_1'$ of $\phi_R(L_1)$. We then have :

$\begin{aligned}
    H_*(L_0) &\simeq LCH^{\ep_0,\ep_0'}_*(L_0, \phi^T_R(L_0)) \\
    & \simeq LCH^{\ep_0,\ep_1'}_*(L_0, \phi^T_R(L_1)) 
    \end{aligned},
    $

Because $ LCH^{\ep_0,\ep_1'}_*(L_0, \phi^T_R(L_1))$ is generated by Reeb chords between $L_0$ and $L_1$ this gives the result.

This proof using contact homology also suggests that there may be no need for a polarisation.

If $(W, \lambda)$ is a Liouville manifold,  $L_1, L_2$ are legendrians in $\mc CW$ we set $$\textbf{Reeb}_a^b(L_1):=\{ (x, z_1, z_2)\in \reeb(L_1, L_1), a<z_2-z_1<b\}.$$The letter $\pi$ without indices will always denote the natural projection $\mc C W\map W$.

If  $g : M \times \R^k\map \R$ is a generating family, denote $$\Sigma_g = \{(x, v) \in M\times \R^k\,|\, \partial_vg = 0\}$$ its critical set, $\iota_g : \Sigma_g\map J^1M$ the associated Legendrian immersion and $\bar \iota_g = \pi\circ\iota_g$.

\section{Preliminary lemmas}

Let $E\cong W \times \R^{N_0}$ be a trivial bundle over $W$ such that $P$  is a polarisation in $\Lambda(TW\oplus E\oplus E^*)$ and $P_{|\lbar}$ is homotopic to $T\lbar\oplus E$ in $\Lambda(TW\oplus E\oplus E^*),$ the Lagrangian Grassmannian bundle of $TW\oplus E\oplus E^*.$

\begin{lemm}
    Let $P_1$ be a polarisation of a symplectic manifold $(W_1, \lambda_1)$, $I_1\subset W_1$ an isotropic submanifold and $L_1 \supset I_1$  a Lagrangian thickening of $I_1$. Assume that $P_{1{|I_1} }$ is homotopic to $T_{|I_1}L_1$ in $\Lambda(W_1)_{|I_1}$. Then there exists a polarisation $P_2$ homotopic to $P_1$ such that $P_{2|I_1} = T_{|I_1}L_1$.
\end{lemm}
\begin{proof}

Let $H_{t}' : I_1\map \Lambda(TW_1)$ be a homotopy between $T_{|I_1}L_1$ and $P_{1|I_1}$. Extend $H_t'$ into a map $H_t : L_1 \map \Lambda(TW_1)$. The homotopy $H$ lifts to a map $\tilde H_t : L_1\map U(TW_1)$ such that $H_t(x)  = \tilde H_t(x)\cdot P_{1,x}$ and $\tilde H_1(x) = \textbf{id}$.

Work on a small Weinstein neighbourhood $V \approx T_{\leq \eta}^*L_1$ of $L_1$. On $V$, the distribution $\mc V = \textbf{ker}(d\pi_{L_1})$ is Lagrangian. Fix an almost complex structure $J$ on $V$ such that $J\mc V_x = T_xL_1$ for $x\in L_1$ and set $\mc H_x = J\mc V_x$ for all $x\in V$. The family of diffeomorphisms $d\pi_{\lbar,x} : \mc H_x\map \mc H_{\pi_\lbar(x)}$ lifts uniquely to a family of symplectomorphisms $\psi_x : T_xV\map T_{\pi_{L_1} (x)}V$ satisfying $\psi_x(\mc V_x) = \mc V_{\pi_{L_1}(x)}$. 
Set
\[
P_{2,x} =
\begin{cases}
P_{1,x} & \text{if } x \notin V, \\[6pt]
(\psi_x)\inv\!\left(\tilde H_{\frac{\|x\|^2}{\eta^2}}(\pi_{L_1}(x))\cdot\big(\psi_x(P_{1,x})\big)\right) & \text{if } x \in V.
\end{cases}
\]
For $x\in L_1$ we have $\psi_x = \textbf{id}$ hence $P_{2,x} = \tilde H_0(x)\cdot(P_{1,x}) = H_0(x) =T_x{L_1}$. Moreover if $\|x\|=\eta$, then $\tilde H_1 = \textbf{id}$ so $P_{2,x}=P_{1,x}$ on $\partial V$.

\end{proof}

Lemma $(3)$ applied to $W_1 = W\times\R^{N_0}\times(R^{N_0})^*$,  $P_1 = P$, $L_1 = \lbar\times \R^{N_0}\times \SingZ$ and $I_1=\lbar$ gives a polarisation $P'$ of $ W\times \R^{N_0}\times(\R^{N_0})^*$ such that for all $x\in \lbar$ : $$P'_x = T_x(\lbar\times \R^{N_0}\times\SingZ).$$

\begin{deff}
    If $f:A\map B$ is a symplectic map between $(A, \omega_A)$ and $(B, \omega_B)$, let $s\nu(f)$ denote the symplectic subbundle of $f^*(TB)$ defined by, for $a\in A :$  $$s\nu(f)_a := \textbf{im}(df_a)^{\perp_{\omega_B}}.$$
\end{deff}
\begin{deff}
    If $(W,\lambda)$ is a Liouville manifold, let 
\[
\begin{aligned}
\iota_\lambda : W &\longrightarrow T^*W \\
x &\longmapsto \lambda_x
\end{aligned}
\]

    denote the graph of $\lambda$.
\end{deff}
\begin{lemm}
    The restrictions $\pi^\perp := d\pi_{W|s\nu(\iota_\lambda)}$ and $\pi^{//} := d\pi_{|\textbf{im}(d\iota_\lambda)} $ induce symplectic vector bundle isomorphisms :     $$(s\nu(\iota_\lambda), -dp\wedge dq) \xrightarrow[]{\pi^\perp}(TW, \omega)$$ and $$(\textbf{im}(d\iota_\lambda), dp\wedge dq) \xrightarrow[]{\pi^{//}}(TW, \omega).$$
\end{lemm}
\begin{proof}
    Set $\iota_\lambda(x) = (x, \lambda_x)$. By the tautological property of $\lambda_{W} = {pdq}$ we have $\iota_{\lambda}^*(pdq) = \lambda$.
    Write $\lambda = \sum_i f_i(q)dq_i$. T hen $d\lambda = \sum_i df_i\wedge dq_i = \sum_{i, j} \dron{q_j}f_idq_j\wedge dq_i$.
If $X_a = \sum_i a_i \dron{q_i}$ is an arbitrary element of $TW$, then $$d\iota_\lambda(X_a) = \sum_ia_i\dron {q_i}+ \sum_{i,j} a_i (\dron{q_i}f_j) \dron{p_j}.$$
One immediately checks that $s\nu(\iota_\lambda)$ is generated by vectors of the form : $$Y_b = \sum_ib_i\dron {q_i}+ \sum_{i,j} b_j (\dron{q_i}f_j) \dron{p_i}.$$
Then 
\begin{align*}
    dp\wedge dq(Y_a, Y_b) &= \sum_{i, j}a_ib_j(\partial_{q_j}f_i-\partial_{q_i}f_j) \\ &= -d\lambda(a\partial_q, b\partial_q)
\end{align*} and $dp\wedge dq (X_a, X_b) = d\lambda(a\partial_q, b\partial_q)$.
We deduce that $d\pi_{W}$ induces symplectic vector bundle isomorphisms between  $(s\nu(\iota_\lambda), dp\wedge dq)$ and $(TW, -d\lambda)$ and between $(\textbf{im}(d\iota_\lambda), dp\wedge dq) $ and $(TW, d\lambda)$.  
\end{proof}

\begin{lemm}
   Let  $(W_1, \lambda)$ be a Liouville manifold equipped with a polarisation $P$. There exists a smooth manifold $M$, an embedding $j : W_1\map T^*M$ that preserves the Liouville form, and a trivial Lagrangian subbundle $\mc I\subset s\nu(j)$ such that for $x \in W_1$,
$$\mc L_x :=\mc I_x\oplus dj(P_x)\xrightarrow[]{d\pi_{M, x}}T_xM$$
is an isomorphism. 
\end{lemm}
\begin{proof}
Let \( p : M \to W_1 \) be a vector bundle and \( z : W_1 \to M \) its zero section. Let $\pi_{W} : T^*W_1 \to W_1$ be the cotangent bundle of \( W_1 \).  The differential \( dp : TM \to TW_1\) induces the dual map
\[
l : (x, \beta)\in T^*W_1 \longrightarrow (z(x), \beta\circ dp)\in T^*M
\]
 over $z$. 
The image of \( l \) consists of covectors based on \( z(W_1) \subset M \) that vanish on $\ker(dp) \subset TM .$ 
 We introduce the following monomorphisms over $l$  :
\[
l_1 : TT^*W_1 \xrightarrow{\, dl \,} TT^*M ,
\]
\[
l_2 : ((x,\beta),(x,d)) \in \pi_{W_1}^*(M) \longmapsto 
\left.\frac{d}{dt}\right|_{t=0} \bigl((x,td),\, \beta\circ dp\bigr) \in TT^*M ,
\]
If $U\subset W_1$is a neighbourhood on which $p\inv(U)\approx U\times \R^n$ is trivial, then
$$l : T^*U\map T^*U\times T^*\R^n$$ is the inclusion map and 
$$l_2: T^*U\times \R^n\map TT^*U\times TT^*\R^n$$
is the product map that sends $T^*U$ on the zero section of $TT^*U$ and $\R^n$ to the zero section of  $T_0T^*\R^n\approx T^*\R^n$.

We deduce that for $x\in T^*U :$  $$s\nu(l)_x = \textbf{im}(l_1^x)^\perp \approx \textbf{im}(l_2^x )\oplus (\textbf{im}(l_2^x ))^* \approx (\pi_{W}^*M)_x\oplus (\pi_{W}^*M)^*_x.$$
 
We know $T{W_1} \simeq P \oplus P^* $. Then  \(s\nu(\iota_{\lambda}) \)  is also isomorphic to  $P \oplus P^* $, by the map $(\pi^{//})\inv\circ \pi^{\perp}$. Fix a vector bundle \( M \) such that \(P \oplus M \) is trivial.
The composite map
\[
j: W \xrightarrow{\iota_{\lambda}} T^*W \xrightarrow {l} T^*M
\]
is then symplectic by lemma $(6)$, and its symplectic normal bundle is
\begin{align*}
s\nu (j) &= s\nu(\iota_\lambda)\oplus \iota_\lambda^*s\nu(l) \\ &= P\oplus P^*\oplus \iota_\lambda^* (\pi_W^*M)\oplus \iota_\lambda^*(\pi_W^*M)^*  \\ &= (P\oplus M)\oplus(P\oplus M)^*
\end{align*}
Set  $\mc I_x = dl\Bigl((\pi^\perp_{W})\inv( P_x')\Bigr) \oplus  \textbf{im}(l_2^x) $, where $P'$ is any polarisation transverse to $P$. In coordinates, $\mc I_x = (\pi^\perp_{U})\inv(P_x')\oplus \R^n\subset T_{p(x)}T^*U\oplus T_0T^*\R^n$  (again $\R^n$ here is the zero section of $T^*\R^n \approx T_0T^*\R^n$) and $dj(P_x) = (\pi^{//})\inv(P_x)\oplus 0$. We deduce that $d\pi_M(\mc I_x\oplus dj(P_x)) = T_{p(x)}U\oplus\R^n$, so $\mc L_x$ is transverse to $\textbf{ker}(d\pi_M)$.
\end{proof}

 We use lemma $(7)$ to find a symplectic embedding $j^N : W\times \R^N\times\R^N\map T^*M$ and a subbundle $\mc I^N \subset s\nu(j^N)$. 
Set $j = j^N_{|W}: W\map T^*M$, $\mc I= \mc I^N\oplus E \subset s\nu(j)$. Then $j$ is symplectic and  $\mc I$ is a trivial lagrangian subbundle of $s\nu (j)$. For $x\in \lbar$, since $P_x'=T_x\lbar\oplus E$, lemma $(7)$ also shows that $T_x\lbar \oplus dj(\mc I)$ is graphical in $T^*M$.

The following lemma follows by Weinstein's neighbourhood theorem (\cite{Weinstein1971}) :
\begin{lemm}
    Let $j : W_1\map  W_2$ be a symplectic embedding, and $\mc I$ a lagrangian subbundle of $s\nu(j)$. Suppose that $\mc I$ is trivial as a bundle over $W_1$.
    There exists a neighbourhood $\mc U$ of $j(W_1)$ in $W_2$, $N\in \mathbb{N}$ and a symplectic embedding $\mc S : \mc U \map W_1\times \R^N\times(\R^N)^*$ such that for all $x\in W_1$, $\mc S_*(\mc I_x) =  \{0\}\times T_0\R^N\times \{0\}$ and  $\mc S(x) = (x, 0, 0)$. We can lift $\mc S$ to a contact embedding $\tilde S : \mc C \mc U \map \mc CW_1\times \R^N\times (\R^N)^*$.
\end{lemm}

Lemma $(8)$ allows us to embed a neighbourhood $\mc U\subset T^*M $ in $ W\times\R^N\times(\R^N)^*$ while sending $j(W)$ to the zero section and $\mc I$ on the tangent spaces of the $\{ x\} \times\R^N\times\SingZ$. 

\section{Construction of the isotopy \texorpdfstring{$\widehat L_s$}{Ls}}

\begin{deff}
    Let $(V, \omega)$ be a symplectic vector space, and $L_1, L_2, L_3\subset V$ three lagrangian subspaces. Suppose that $L_1\pitchfork L_3$ and $L_2\pitchfork L_3$. Then $L_1$ is the graph of a symmetric linear map $s : L_2\map L_3\approx_\omega L_2^*$.
    Set $$\textbf{Ind}(L_1, L_2, L_3) := n_+-n_-$$ where $n_-$ (resp. $n_+$) is the number of negative (resp. positive) eigenvalues of  $s$ in any basis.
\end{deff}

We denote $\textbf{Red}(L_1, I)$ the symplectic linear reduction of $L_1$ with respect to an isotropic subspace $I$.

\begin{lemm}    Let $I$ be an exact isotropic submanifold of a cotangent bundle $T^*M$. Let $\tilde I$ denote a lift of $I$ in $\mc CM.$ Let $H\subset T^*M$ be a coisotropic submanifold, $L$ a lagrangian thickening of $I$ and $K\subset T^*M$ a compact subspace whose interior contains $I$. Suppose that $L\pitchfork H$ with $L\cap H = I$ and that $L$ is graphical in $T^*M$.
    There exists a generating family $g: M\times \R^k\map \R$ of a submanifold $\widehat L=\widehat L^+\cup \widehat L^- \subset J^1M$ satisfying the following : 
\begin{enumerate}
  \item  There exists $\mu>0$ such that $\reeb_0^\infty(\Lch) = \{(x, z-\mu, z+\mu) ; (x, z) \in \tilde I\}$.
  \item $\widehat{L} \subset K\times \R$; 
  \item $\pi(\Lch^\pm) \cap H = I$ ;
  \item $\pi(\Lch^\pm) \pitchfork H$ at $I$.
  \item For $p\in I\subset T^*M$  : $$\textbf{Ind}(\textbf{Red}(T_p\pi(\Lch^\pm), T_pI), \textbf{Red}(T_pL, T_p I), \textbf{Red}((T_pH)^\perp\oplus T_pI, T_pI)) = \pm(\textbf{dim}(M)-\textbf{dim}(I)).$$
  
\end{enumerate}
\end{lemm}
\begin{proof}
We start by constructing auxiliary functions $f$ and $A$ that will help in defining $g$.

After shrinking $L$, we can assume that a contact lift $\tilde L$ of $L$ deformation retracts onto $\tilde I$ and that $L\subset K$. Since $L$ is graphical we can find a function $f: M\map \R$ whose $1-$jet contains $\tilde L$.

Choose a function $A \colon M \to \mathbb{R}$ such that
\begin{itemize}
    \item $A$ is of compact support ;
    \item $\displaystyle
        \pi(I)
        = \left\{
            x \in \operatorname{Crit}(A)
            \,\middle|\,
            A(x) > \frac{1}{2}
          \right\}$ ;
    \item $A(\pi(I)) = \{2\}$ and $A$ is Morse-Bott at $A\inv(\{2\})$.
\end{itemize}

Define the following generating function (with only one external variable $w$):
\[
g(x;w)=f(x)+w^3+3\ep w-3\varepsilon wA(x),
\]
where \(\varepsilon>0\) will be chosen later.

 We have
\begin{itemize}
    \item $\Sigma_g=\{(x,w)\mid w^2=\varepsilon (A(x)-1)\},$
    \item $\mathbf{Crit}(\partial_w g)=\textbf{Crit}(A)\times\SingZ,$
\end{itemize}
so that
\[
\mathbf{Crit}(\partial_w g)\cap \Sigma_g
=\textbf{Crit}(A)\times\{0\}\cap \{(x,w)\mid w^2=\varepsilon (A(x)-1)\}
= \varnothing,
\]
since \(A(\pi_M(I))=\{2\}\). Thus, \(g\) is really a generating family.

Let $\widehat L = \iota_g(\Sigma_g)$  denote the legendrian generated by $g$ and $$\widehat L^{\pm} := \iota_g(\Sigma_g\cap M\times\R^\pm)$$ the subsets of $\Lch$ corresponding to  $\pm w\geq 0$. 

Let $\bar\iota_{g}^\pm : x\in A\inv([1, 2])\mapsto \bar \iota_g(x, \pm\sqrt{\ep(A(x)-1)})$ be a parametrisation of $\pi(\Lch^\pm)$ over $x\in M$. For $(x, w)\in \Sigma_g$ and $\mp w\geq 0$ : $$\iota_{g}(x, w) = (\bar\iota_g^\pm(x), f(x)\mp2\ep^{\frac32}(A(x)-1)^{\frac32}).$$
We deduce that $\iota_{g|\Sigma_g}$ is injective and the reeb chords in $\reeb(\Lp, \Lm)$ are of the form $((x, df(x)), f(x)-2\ep^{\frac32}, f(x)+2\ep^{\frac{3}{2}})$ for $x\in I$ from which $(1)$ follows for $\mu = 2\ep^\frac32$.
Moreover, on the set $\{A\geq 1\}$ : $$\bar \iota_{g}^\pm\xrightarrow[]{\ep\map 0} \bar\iota_{f}$$ and on the set $\{A > 1\} : $ $$\bar\iota_{g}^\pm\xrightarrow[\mc C^1]{\ep\map 0} \bar\iota_{f}.$$
Since $\textbf{im}(\bar\iota_f)=L\subset K $ above the set $\{A > 1\},$ we deduce $(2)$ is true for $\ep$ sufficiently small. Because $L$ intersects $H$ transversally in $I$,  $(3)$ and $(4)$ follow for $\ep$ sufficiently small.

We show $(5)$. Let $p\in  I$, and find coordinates $(s, t)$ around $\pi_M(p)$ in $M$ for which  $A(s, t) = 2-\|t\|^2. $ Then $T_p\pi(\Lp)/T_pI=\textbf{Vect}_i(\partial_{t_i}\bar\iota_g^+(s, 0))/T_pI$. 
Furthermore
\begin{align*}
    \partial_{t_i}\bar\iota_g^+(s, 0) &=  \partial_{t_i}\bar\iota_f(s, 0)-\partial_{t_i}|_{t=0}\bigl( 0,0, 0, 3\ep \sqrt{\ep (A(s, 0)-1)} \partial_tA(t)\bigr)\\
    &=  \partial_{t_i}\bar\iota_f(s, 0)+\partial_{t_i}|_{t=0}\bigl( 0, 0, 0, 6\ep^{\frac32} t \bigr)\\
    &\in T_pL/T_pI\oplus \Bigl(\textbf{ker}(d\pi_{M,p})\cap (T_pI)^{\perp}\Bigr)
\end{align*}
with $\omega(\partial_{t_i}\bar\iota_f(s, 0),\partial_{t_i}|_{t=0}\bigl( 0, 0, 0, 6\ep^{\frac32} t \bigr))=6\varepsilon^{\frac32}>0$.

In the end this shows $$\textbf{Ind}(\textbf{Red}(T_p\pi(\Lp), T_pI), \textbf{Red}(T_pL, T_p I), \textbf{Red}(\textbf{ker}(d\pi_{M,p}), T_pI)) = \textbf{dim}(M)-\textbf{dim}(I).$$
Since $T_p\Lp$ is $\ep-$close to $T_pL$, the index does not depend on the third lagrangian for $\ep$ sufficiently small. This proves $(5)$, the calculation for $\Lm$ being essentially the same.

\end{proof}

Choose $W_c\subset W$ a Liouville domain of $W$ containing $\lbar$ and $R\in \R^+$ such that $W_c\times D_R\times D_R\subset \mc S(\mc U)$.
Apply lemma $(10)$ with $\tilde I = L_0$, $L = \mc{S}\inv(\lbar\times \R^N\times\{0\})$, $H = \SSS\inv(W\times \SingZ\times (\R^N)^*)$ and $K = \mc S\inv(W_c\times D_R \times D_R)$, to obtain a generating family $g$ of a legendrian $\Lch = \Lp\cup\Lm.$ The hypothesis of $L$ being graphical near $I$ is satisfied because by definition, for $x\in \lbar$, the tangent space $T_{(x, 0, 0)}(\lbar\times \R^N\times\SingZ)$ is $d\mc S(\mc L_x)$.
Set $\widehat K = \tilde{\mc S} (\widehat L)$ and $\Kch^\pm = \tilde{\mc S} (\widehat L^\pm)$.

\begin{figure}[ht]
    \centering
    \includegraphics[width=0.8\textwidth]{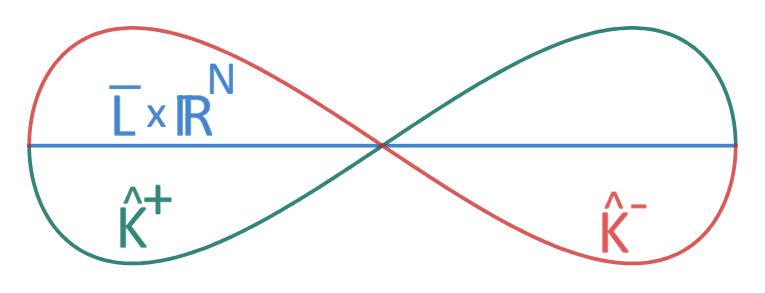}
    \caption{Tranche de $\pi(\Kch)$ dans $\{x\}\times\R^N\times(\R^N)^*$}
    \label{fig:geogebra2}
\end{figure}

\begin{lemm}
There exists a positive real number $\eta$ and a parametrisation
\[
\begin{aligned}
s^+ : \lbar\times D^N_{_\leq \eta}  &\map W\times\R^N\times (\R^N)^* \\
(p, t) &\mapsto (s_1^+(p, t), t, S^+(p)t+R^+(p, t)t)
\end{aligned}
\]
of a neighbourhood of $\lbar$ in $\pi(\Kp)$ such that for all $(p, t)\in \lbar\times D^N_{\leq\eta}:$
\begin{itemize}
\item $s_1^+(p, 0) = p$
\item $S^+ : \lbar \map S_N^{++}$
\item $R^+(p, 0) = 0\in \mc M_N(\R)$
\item $|t^TR^+(p, t)t|\leq \frac {t^TS^+(p)t}2$.
\end{itemize}

\end{lemm}
\begin{proof}
    
Let $U_L\approx T^*_{<v}L \subset W$ denote a small Weinstein neighbourhood of $\lbar$ in $W$.
Above  $U_L$ in $W\times \R^N\times (\R^N)^*$ there is a projection (induced by $T_{\leq v}^*\lbar\map \lbar$) : $$\textbf{pr} : U_L\times \R^N\times (\R^N)^*\map \lbar\times\R^N.$$
For $p\in L$, the restriction of  $\textbf{pr}$  to $T_p \pi(\Kp)$ is nonsingular because $\pi(\Kp)$ intersects $ W\times  \{0\}\times (\R^N)^*$ transversally and $T_p\pi(\Kp) \supset T_p\lbar$ by the points $(3)$ and $(4)$ of lemma $(10)$.

We can thus parametrise a neighbourhood of  $\lbar$ in $\pi(\widehat K ^+)$ by a section of $\textbf{pr}$ : $$(s^+) : L\times D_{\leq \eta}\map W \times \R^N\times (\R^N)^*.$$
Since $\pi(\Kp)\cap (W\times\{0\}\times(\R^N)^*) = L\times\{0\}\times\{0\}$ we have $\pi(\widehat K^+ )\cap (W\times D_{\leq\eta}^N\times(\R^N)^*) = \textbf{im}(s^+)$ for $\eta$ sufficiently small.

Write $s^+ = (s_1^+, s_2^+, s_3^+)$. By construction, for $p\in \lbar$, $\textbf{pr}\circ s^+(p, t) = (p, t)$ so $s_2^+(p, t) = t$. Since $\lbar\subset \pi(\Kp)$we have moreover $s^+(p, 0) = (p, 0, 0)$ and thus $$ds^+_{(p, 0)}(T_pL\times \SingZ) = T_pL\times\SingZ\times\SingZ.$$
and \begin{align*}
    \partial_ts^+(p, 0) &\subset (T_pL\times\SingZ\times\SingZ)^{\perp\omega} \\ &=T_pL\times T_0\R^N\times T_0(\R^N)^*.
    \end{align*}
This shows $\partial_ts_1^+(p, 0)\in T_pL$ and since $\textbf{pr}\circ s_+ = \textbf{id}$, we also have $\partial_ts^+_1 \in \textbf{ker}(d\pi_{L, p})$. Therefore $\partial_t s_1^+(p, 0) = 0$.

Let $S^+(p) := \partial_ts_3^+(p, 0)\in \mc M_N(\R)$ denote the jacobian matrix of $s_3^+(p, -): D_{\leq\eta}\map\R^N$. We can write $\partial_ts^+(p, 0) = (0, \partial_t, S^+(p)\partial_\tau)$. Because the vector space  $ds_{(p, 0)}^+(\SingZ\times T_0\R^N)=\SingZ \times \textbf{graph}(S^+(p))$ is isotropic, $S^+(p)$ is a symmetric map.

We now show $S^+(p)\in S_N^{++}$ . Identify $T_pL^\perp /T_pL$ with $\R^N\times (\R^N)^*.$ Then :
\[
\begin{aligned}
\textbf{red}(T_p\pi(\Kp), T_pL)
&=\textbf{red}(\textbf{im}(ds_{(p, 0)}^+), T_pL)  \\
&=\textbf{red}(T_pL\times \textbf{graph}(S^+(p)), T_pL)  \\
&\simeq \textbf{graph}(S^+(p)).
\end{aligned}
\]
Furthermore
$$\textbf{red}(T_p(L\times\R^N\times\SingZ), T_pL) = T_0\R^N\times\SingZ.$$
and $$\textbf{red}(T_p(\{p\}\times\SingZ\times(\R^N)^*), T_pL) = \SingZ\times T_0(\R^N)^*$$
Conclusion $(5)$ of lemma $(10)$ thus gives that $S^+(p)\in S_N^{++}$.

Since $s_3^+(p, 0) = 0$, the Hadamard lemma allows us to write $s_3^+(p, t) = S^+(p)t+R^+(p, t)t$ for some function $R^+ : L\times D_{\leq\eta}\map \mc M_N(\R)$. Because $\partial_ts_3^+(p, 0) = S^+(p)\partial_\tau$ we have that $R^+(p, 0) =0$.

A Taylor expansion gives
\[
|t^TR^+(p,t)t|
\leq C'\|t\|^3
\]
for some $C'>0$. Moreover $t^TS^+(p)t\geq C\nrmtc$ for some $C>0,$ so that we can conclude $$|t^TR^+(p, t)t|\leq \frac {t^TS^+(p)t}2$$ by taking $\eta\leq \frac{C}{2C'}$.

\end{proof}
Let $z^+ : \lbar\times D^N_\eta\map \R$ denote the function such that $(s^+, z^+)$ is a local parametrisation of $\Kp$.
We similarly construct a parametrisation $s^-$ of a neighbourhood of $\lbar$ in $\pi(\Km)$, and analogous functions $s^-_1(p, t), R^-(p, t)$,  $S^- : L\map S_n^{--}$ and $z^- : \lbar\times D^N_\eta\map \R$.

\begin{lemm}
        We construct a Legendrian isotopy $\widehat K_s$ equipped with a function $ n' : [0, 1]\map \R^*_+$ such that :
    \begin{enumerate}
        \item $\widehat K_0 = \widehat K$
        \item For $s\in [0, 1]$, $\Kch_{s}$ does not have Reeb chords of length $n'_s$
        \item $\#\reeb_{n'_1}^\infty(\Kch_1) = \#\reeb(L_0, L_1)$
        \item  Reeb chords in $\reeb_{n'_1}^\infty(\Kch_1)$ are non degenerate.
        \item $\reeb_{n_0'}^\infty(\widehat K_0)$ is diffeomorphic to $L_0$.
    \end{enumerate}

\end{lemm}

\begin{proof}
    
Start by choosing a family of contactomorphisms $(\varphi_t)_{t\in [0, 1]} : \mc CW\map \mc CW$ satisfying $\varphi_0 = \textbf{id}_{\mc CW}$ and $L_t = \varphi_t(L_0)$. Let $m_t :\mc CW\map \R^*_+$ be the function such that $\fee_t^*(dz-\lambda) = m_t\cdot(dz-\lambda)$. We consider the following families of contactomorphisms on $(W\times \R^N\times(\R^N)^*\times\R, dz-\lambda-\tau dt)$:
\[
\psi_s(x,t,\tau,z)
= \bigl(x,\; t,\; \tau+2Ms\,\gamma'(\|t\|^2)t,\; z+Ms\,\gamma(\|t\|^2)\bigr),
\]
where $M\in\mathbb{R}^+$ and $\gamma:\mathbb{R}^+\to\mathbb{R}$ is a decreasing function such that
\[
\gamma(s)=0 \quad \text{if } \sqrt{s}>\tfrac{\eta}{2},
\qquad
\gamma(s)=1-s \quad \text{if } \sqrt{s}<\tfrac{\eta}{3},
\]
and
\[
\chi_s(x,t,\tau,z)
= \bigl(\feetild^x_{s r(\|t\|^2)}(x,z),\; t,\;
m_{sr(\|t\|^2)}\tau+(dz-\lambda)(\partial_t[\feetild_{s r(\|t\|^2)}(x, z)]),\;
\feetild^z_{sr(\|t\|^2)}(x,z)\bigr),
\]
where $\feetild_s(x, z+\mu) = (\fee_s^x(x, z), \fee_s^z(x, z ) +\mu),$ $\mu$ was defined in lemma $(10),$  and $r:\mathbb{R}^+\to\mathbb{R}$ satisfies 
\[
r(s)=0 \text{ if } \sqrt{s}>\tfrac{\eta}{4},
\qquad r(0)=1.
\]
  
Since $\chi^\tau_s-m_{sr(\|t\|^2)}\tau$ vanishes at $t=0$, Hadamard's lemma provides a function :
\[
\Gamma_s : W\times\mathbb{R}^N\times\R \longrightarrow \mathcal M_N(\mathbb{R})
\]
so that :
\[
\chi_s(x,t,\tau,z)
= \bigl(\feetild^x_{s r(\|t\|^2)}(x,z),\; t,\;
m_{sr(\|t\|^2)}\tau+\Gamma_s(x,t,z)t,\;
\feetild^z_{sr(\|t\|^2)}(x, z)\bigr).
\]
By construction, for $\|t\|>\tfrac{\eta}{4}$, we have $\Gamma_s(x,t,z)=0$.

Define an isotopy of Legendrians $\bigl(\widehat{K}_s\bigr)_{s\in[0,1]}$ by :

\[
\widehat{K}^+_s := \Kp,
\]

\[
\Km_s :=
\begin{cases}
\psi_{2s}\!\bigl(\Km\bigr), 
& \text{if } s\in\left[0,\tfrac12\right],\\[0.3em]
\psi_1\circ\chi_{2s-1}\!\bigl(\Km\bigr),
& \text{if } s\in\left[\tfrac12,1\right],
\end{cases}
\]

$$\Kch_s := \Kp_s\cup \Km_s.$$

\begin{figure}[ht]
\centering
\begin{minipage}{0.48\textwidth}
    \vspace{2.13cm}
    \centering
    \includegraphics[width=\linewidth]{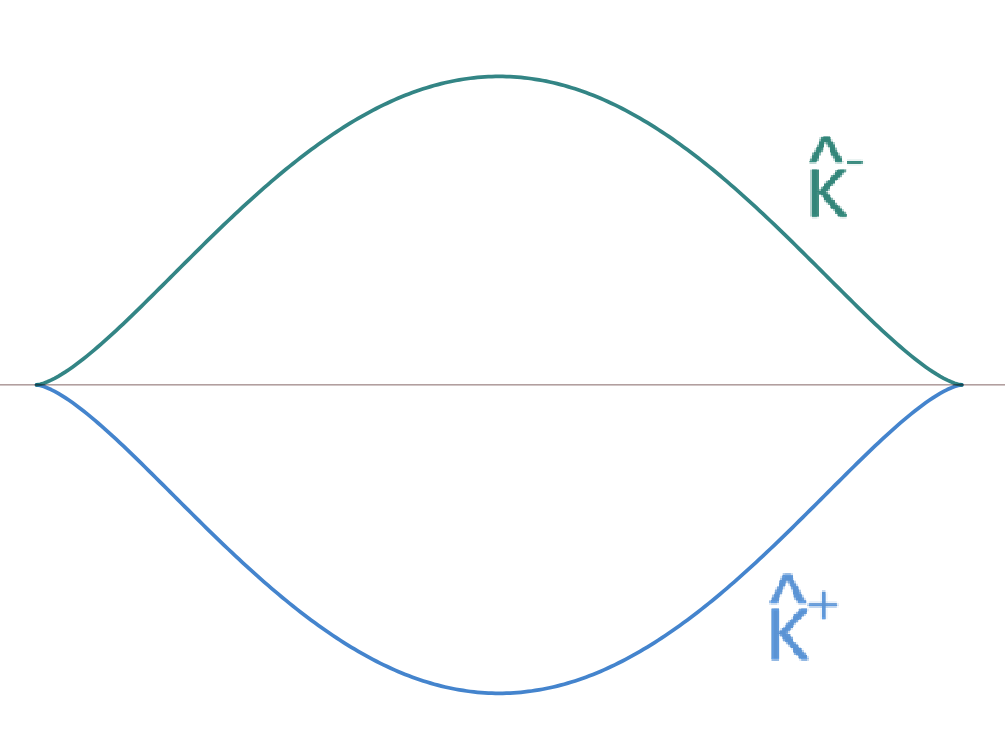}
    \caption{Slice of the front of $\Kch$.}
\end{minipage}
\hfill
\begin{minipage}{0.48\textwidth}
    \centering
    \includegraphics[width=\linewidth]{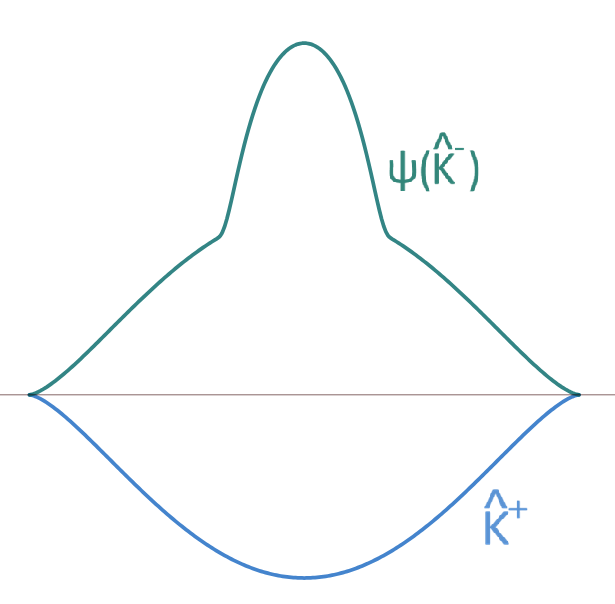}
    \caption{Slice of the front of $\Kch_{\frac12}$.}
\end{minipage}

\end{figure}

\begin{figure}[ht]
    \centering
    \includegraphics[width=0.4\textwidth]{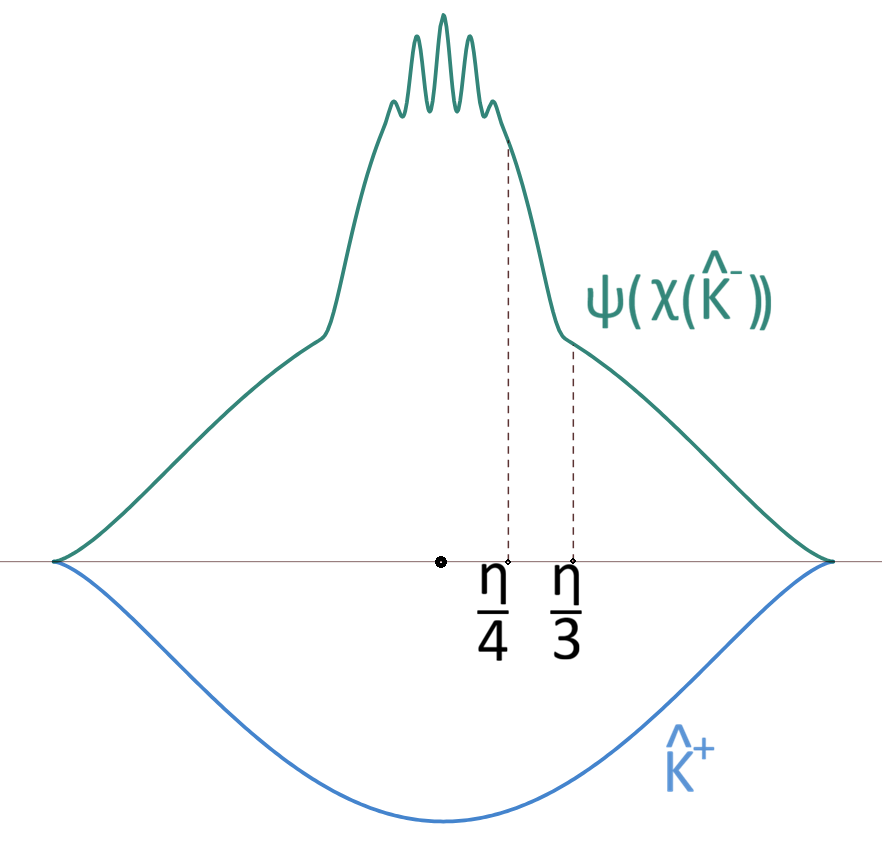}
    \caption{Slice of the front of $\Kch_1$. The oscillation represents $\chi$ and is also happening transversely to this slice.}
    \label{fig:geogebra}
\end{figure}

We now show that the isotopy $\widehat K_s$ has the required properties.

We decompose
$$\reeb_0^\infty(\Kch_s) = \reeb_0^\infty(\Kp_s) \cup\reeb_0^\infty(\Km_s) \cup\reeb'(\Kp_s, \Km_s),$$
where we write $ \reeb'(\Kp_s, \Km_s):=\reeb_0^\infty(\Km_s, \Kp_s)\cup \reeb_0^\infty(\Kp_s, \Km_s)$ for notational convenience.
We have that \begin{equation}
\textbf{Reeb}_0^\infty(\Kp_s) = \textbf{Reeb}_0^\infty(\Kp_0) = \emptyset.
\label{eq:reeb1}
\end{equation}
Moreover, since $\psi_s$ is the contact lift of a symplectomorphism, we have for $s\in [0, \frac12]$ :
\begin{equation}
    \textbf{Reeb}_0^\infty(\Km_s) = \textbf{Reeb}_0^\infty(\psi_{2s}(\Km_0)) = \reeb_0^\infty(\Km_0)=\emptyset.
    \label{eq:reeb2}
\end{equation}
Take now $s\in [\frac12, 1]$ and set $s' = 2s-1 \in [0, 1]$. If $(x, z_1, z_2) \in \reeb_0^\infty(\Km_s)$, then $(x, z_i) = \psi_1\circ\chi_{s'} (x_i, t, \tau_i, z_i')$ for $(x_i, t, \tau_i, z_i')\in \Km_0, (i=1, 2)$. By the definitions of $\psi$ and $\chi$, $z_i = \feetild^z_{s'r(\|t\|^2)}(x_i, z_i')+M\gamma(\|t\|^2)$ and so : 
\begin{align*}
    |z_1-z_2| &= |\feetild^z_{s'r(\|t\|^2)}(x_1, z_1')-\feetild^z_{s'r(\|t\|^2)}(x_2, z_2')| \\ &< 2\displaystyle\sup_{s\in [0, 1],(x, t, \tau, z)\in \Kch_0}(|\feetild^z_s(x, z)|) =:B.
\end{align*}
In the end 
\begin{equation}
    \reeb_{B}^\infty(\Km_s) = \emptyset.
    \label{eq:reeb3}
\end{equation}

Take $s\in [0, \frac12]$ and $((x, t, \tau), z_1,z_2) \in \reeb'(\Kp_s, \Km_s).$ Write $((x, t, \tau), z_2) = \psi_{2s}(x, t, \tau_2, z_3)
$ with $(x, t, \tau_2, z_3)\in \Km_0$. Then by the definition of $\psi$, $$\tau =\tau_2 +4sM\gamma'(\nrmtc)t,$$
$$z_2 = z_3+2Ms\gamma(\nrmtc).$$
If $\tau\neq \tau_2$, then $\nrmt<\frac\eta2$ which implies $(x, t, \tau)\in \textbf{im}(s^+) $ and $(x, t, \tau_2) \in \textbf{im}(s^-)$. Thus for  $p, p'\in \lbar$, $$\tau_2-\tau  = S^-(p)t+R^-(p, t)t -S^+(p')t-R^+(p', t)t.$$ 
These equalities imply : 

\[
\begin{aligned}
0&\leq -4sM\gamma'(\nrmtc)\nrmtc \\
&= t^T(\tau_2-\tau)\\
&= t^T(S^-(p) + R^-(p,t) - S^+(p') - R^+(p',t))t \\
&= t^T(S^-(p) +R^-(p,t))t - t^T(S^+ (p')+R^+(p',t))t \\
&\leq t^T(\frac12S^-(p) -\frac12 S^+(p'))t \\
&\leq 0,
\end{aligned}
\]

The last inequality is strict for $t\neq 0$, so $t=0$ and $\tau=\tau_2$. It follows that $((x, t, \tau), z_1, z_3)\in \reeb'(\Kp_0, \Km_0)$and $z_2-z_1 = 2Ms\gamma(\nrmtc) +z_3-z_1\geq 2Ms+2\mu$ (using lemma $(10.1)$). In the end, for $s\in [0, \frac12] :$
\begin{equation}
    \reeb_{-\infty}^{\mu+2Ms}(\Kp_{s}, \Km_{s})= \emptyset.
    \label{eq:reeb4}
    \end{equation}

Let $s\in [\frac12, 1]$, $s' = 2s-1$ and $\bigl((x, t, \tau), z, z_2\bigr) \in \reeb'(\Kp_s, \Km_s)$. Let $(x', t, \tau', z')\in \Km_0$ such that \begin{equation}(x, t, \tau, z_2) = \psi_1\circ\chi_{s'}(x', t, \tau', z').
\label{eq:reebeq}
\end{equation} 
If $\|t\|>\tfrac{\eta}{2}$, then $r(\nrmtc) = \gamma(\nrmtc) = 0$ and so  $\psi_1\circ \chi_{s'}$ fixes $(x', t, \tau', z')$. Since moreover $\Km_0 = \Km_s$, we have that $\bigl((x, t, \tau), z, z_2\bigr) \in \reeb'(\Kp_0, \Km_0)$ and so $t=0$ which is absurd.

If $\|t\|<\eta$, then $(x, t, \tau, z)$ is in the image of $(s^+, z^+)$ and $(x', t, \tau', z')$ is in the image of $(s^-, z^-).$ Hence there exist $p,p'\in \lbar$ with $z' = z^-(p', t)$, $z = z^+(p, t)$ and, by \eqref{eq:reebeq}
\[
\psi_1\bigl(\chi_{s'}(s^-(p',t), z')\bigr)=(s^+(p,t), z).
\]
Identifying the $x$ and $\tau$ coordinates, this equates to
\[
\fee^x_{s'r(\|t\|^2)}\bigl(s_1^-(p',t), z'\bigr)=s_1^+(p,t),
\tag{$E_1$}
\]
and
\begin{align*}
\bigl(m_{s'r(\|t\|^2)}S^-(p')+m_{s'r(\|t\|^2)}R^-(p',t)+\Gamma_{s'}(p',t,z')\\+2M\gamma'(\|t\|^2)-S^+(p)-R^+(p,t)\bigr)t=0.
\tag{$E_2$}
\end{align*}
Taking $M> \sup_{s, p, t, z} \|\Gamma_s(p, t, z)\|_2$, since $\Gamma \neq 0 \implies \gamma' = -1$ we have that $$-2M\gamma'\geq\|\Gamma\|_2.$$   Multiplying $(E_2)$ on the left by $t^T$, we get

\begin{align*}
0 ={}&t^T\bigl(m_{s'r(\|t\|^2)}S^-(p')+m_{s'r(\|t\|^2)}R^-(p',t)+2M\gamma'(\|t\|^2) \\&-S^+(p)-R^+(p,t)+\Gamma_{s'}(p', t, z')\bigr)t\\
={}&t^T\bigl(m_{s'r(\|t\|^2)}(S^-(p')+R^-(p',t))\bigr)t-t^T(S^+(p)+R^+(p,t)))t \\&+2M\gamma'(\|t\|^2)\nrmtc+t^T\Gamma_{s'}(p', t, z')t\\
\leq{}& t^T\bigl(\frac12m_{s'r(\|t\|^2)} S^-(p')-\frac12 S^+(p)\bigr)t\\&+2M\gamma'(\|t\|^2)\nrmtc+\|\Gamma_{s'}(p', t, z')\|_2\nrmtc\\
\leq{}& \frac12 \min(1, m_{s'r(\|t\|^2)}) t^T( S^-(p')- S^+(p))t.\\
\leq {}&0.
\end{align*}

The last inequality is strict for $t\neq 0$, so $t=0$ and $(E_1)$ equates to $\feetild_{s'}^x(p',  z') = p.$ Since $(p', 0, \tau', z') = (x', 0, \tau', z') \in \Km$, we have $(p', z'-\mu)\in L_0$, and so $$(p, z+\mu,\fee^z_{s'}(p', z'-\mu))\in \reeb(L_0, L_{s'}),$$ with 
\begin{align*}
    z_2-z &= z\bigl(\psi_1\circ\chi_{s'}(x', 0, \tau',z')\bigr)-z
\\ &= \feetild_{s'}^z(x', z')+M-z \\& >M-2B.
\end{align*}
Hence \begin{equation}
    ((x, t, \tau), z_2, z) \in \reeb^\infty_{M-2B}(\Km_{s'}, \Kp_{s'})
    \label{eq:reeb5}
\end{equation}
and $\#\reeb_{M-2B}^\infty(\Km_{1}, \Kp_{1}) = \# \reeb(L_0, L_{1}) $. 
In the end, taking $M > 3B$, and using \eqref{eq:reeb1}, \eqref{eq:reeb2}, \eqref{eq:reeb3}, \eqref{eq:reeb4} and \eqref{eq:reeb5} we can choose a regular level $n'_s$ with $n'_s = \mu+2Ms$ for $s<\frac{1}{2}$ and $n'_s  =M-2B$ for $s> \frac12+\delta$ for some small $\delta$. The isotopy $\Kch_s$ satisfies $(1), (2), (3)$ and $(5).$

We now show $(4).$ The Reeb chords in $\reeb_{n'_1}^\infty(\Kch_1)$ are of the form $(p, p',  0)\in L\times L\times \R^N$ for $\varphi_1(p) = p'$. Around this point, the Reeb chords are defined by the equations $(E_1)$ and $(E_2)$.
We seek a tangent vector $X$ common to $\pi(\Km_1)$ and $\pi(\Kp_1)$ at $(p, p',  0)$. 
The equation $(E_2), $ of the form $A(t)t=0$ with $A(0)$ nonsingular, shows that $X$ must be tangent to $t=0$. When $t=0$  the equation $(E_1)$ reduces to $\fee^x_{1}\bigl(p', z'\bigr)=p.$ A common tangent vector to $\pi(\Km_1)$ and $\pi(\Kp_1)$ would thus be tangent to  $\pi(L_0)$ and $\pi(L_1)$, which implies $X=0.$ Therefore the Reeb chords of $\reeb_{n'_1}^\infty(\Kch_1)$ are non degenerate.
\end{proof}

\begin{lemm}
There is a continuous function  $n : [0, 1] \map \R^*_+$ and an isotopy $\Lch_s$ such that :
\begin{enumerate}
        \item $\widehat L_0 = \widehat L$
        \item $\Lch_s$ does  not have Reeb chords of length $n_s$.
        \item $\#\reeb^\infty_{n_1}(\Lch_1) = \#\reeb(L_0, L_1)$.
        \item  Reeb chords in $\reeb_{n_1}^\infty(\Lch_1)$ are non degenerate.
        \item $\reeb_{n_0}^\infty(\widehat L_0)$ is diffeomorphic to $L_0$.
\end{enumerate}
\end{lemm}
\begin{proof}
    
Let $\mathcal{F}_\beta : W\times \R^N\times (\R^N)^*\times \R\map W\times \R^N \times (\R^N)^*\times\R$ be the Liouville flow. Then $\mc F_\beta^*\lambda = e^{\beta}\lambda$, $\mc F_\beta^*(\tau dt) = e^\beta\tau dt$ and  $\mc F_\beta(\mathcal{U}) \supset \mc F_\beta(W_c)\times D_R\times  D_{e^\beta R}$. Since $\Kch_0\subset W_c\times D_R\times D_R$ and the isotopy $\Kch_s$ preserves the $t$ coordinate, we can choose $\beta$ large enough so that $\mc F_\beta(\mathcal{U})$ contains the whole isotopy $\Kch_s $. Define
\[
\widehat{L}_s :=
\begin{cases}
\mc S\inv ( \mc F_{-2s\beta}\bigl(\widehat{K}\bigr)), 
& \text{if } s\leq \frac12,\\[0.3em]
\mc S\inv (\mc F_{-\beta}(\widehat K_{2s-1})),
& \text{if } s\geq \frac12.
\end{cases}
\]
We have that $\mathcal F_{-\beta}( x)\in \mathcal{U} \iff x\in \mathcal{F}_\beta(\mathcal{U})$ so $\widehat{L}_s$ is well defined. 
The function
\[
n_s =
\begin{cases}
e^{-2s\beta}n'_s & \text{if } s \leq \frac12,\\[4pt]
e^{-\beta}n'_s & \text{if } s \geq \frac12,
\end{cases}
\]
then satisfies the desired properties.

\end{proof}

\section{Generating family homology of \texorpdfstring{$\widehat L_s$}{Ls}}
 
Proposition $3.1.5$ in \cite{Eliashberg1998-mz} shows that the isotopy $(\Lch_s)_{s\in [0, 1]}$ is generated by an isotopy of generating families  $(g_s)_{s\in[0,1]}$ of the form $g_s(x, w, v) = f(x)+w^3+3\ep w-3\ep w A(x) + Q(v)+\rho_s(x, v, w)$ where $Q$ is a non degenerate quadratic form, $v\in \R^k$, $\rho_s$ is of compact support and $\rho_0  = 0$.

\begin{prop}
    There exists an isotopy of generating families $h_s(x, w, v) = f(x)+w^3+3w\ep +Q(v) +\rho_s(x, v, w)+3 w\eta(x, v)$ with $\rho_0 = 0$ and $\eta$ of compact support that generates $(\Lch_s)_{s\in [0, 1]}$. 
\end{prop}
\begin{proof}
Choose a cutoff function $\chi: \R \map [0, 1]$ such that  $\chi_{|]-\infty, 1]} = 1$ and $\chi_{|[2, \infty]} = 0.$ Define for $v\in \R^k$ $$\chi_R(v) = \chi(\|v\|-R),$$ so that $\|\nabla \chi_R\|_\infty$ is bounded uniformly in $R$.
Define the isotopy 
$$
g_{s,u}(x,v,w)
=
f(x)+w^3+3\varepsilon w
-3\varepsilon w A(x)\bigl((1-u)+u\chi_R(v)\bigr)
+Q(v)+\rho_s(x,v,w)
.$$ We show that for $R$ sufficiently large, $h_s := g_{s, 1}$ is a generating function equivalent to $g_s$, that therefore also generates $\Lch_s$.

    We start by checking that the singular space $\Sigma_{s, u}$ of $g_{s, u}$ is independent of $u$ for $R$ large enough. For $(x, v, w)\in M\times\R^k\times\R$, the set $\Sigma_{s, u}$ is cut out by the equations

\[
\left\{
\begin{aligned}
(E_v)\qquad
\partial_v g_{s,u}(x,v,w)
&=
-3\varepsilon w\,uA(x)\,\partial_v\chi_R(v)
+\partial_v Q(v)
+\partial_v\rho_s(x,v,w),
\\[0.5em]
(E_w)\qquad
\partial_w g_{s,u}(x,v,w)
&=
3w^2+3\varepsilon
-3\varepsilon A(x)\Bigl(
1-u+u\chi_R(v)
\Bigr)
+\partial_w\rho_s(x,v,w).
\end{aligned}
\right.
\]

For $\|v\| < R$, since $g_{s, u}(x, v, w) = g_s(x, v, w)$ we have $(x, v, w)\in \Sigma_{s, 0}$. By taking $R$  large enough so $\Sigma_{s, u}\subset \{\|v\|<R\}$ we can ensure $$\Sigma_{s, u} \cap \{\|v\| < R\} = \Sigma_{s, 0}.$$
    For  $\|v\| \geq R$, the equation $(E_w)$ gives $$3w^2 \leq 3\ep +3\ep \|A\|_\infty+\|\nabla\rho\|_\infty, $$and $(E_v)$ yields, using that $Q$ is non degenerate :
    
    \begin{align*}
      R &\leq  \|v\| \\ &\leq
   C \|\partial_vQ(v)\| \\ &\leq 3C|w|\|A\|_\infty\|\nabla\chi_R\|_\infty+C\|\nabla\rho\|_\infty \\ &\leq 3C\sqrt{ 3\ep +3\ep \|A\|_\infty+\|\nabla\rho\|_\infty}\|A\|_\infty\|\nabla\chi_R\|_\infty+C\|\nabla\rho\|_\infty
    \end{align*}
    
    for some constant $C>0$. Therefore $\Sigma_{s, u} \cap \{\|v\| \geq R\} = \emptyset$ for $R$ sufficiently large and so $$\Sigma_{s, u} = \Sigma_{s, 0}.$$

It remains to be shown that $(g_{s, u})_{u\in [0, 1]}$ can be realised by a fibered diffeomorphism.
We set \[X_{s, u}=
\begin{cases}
0, & \|v\|<R,\\[2mm]
-\displaystyle\frac{\partial_u g_{s, u}}
{\|\nabla_{w,v}g_{s, u}\|^2}\,\nabla_{w,v}g_{s, u},
& (x,w,v)\notin\Sigma_{s, u},
\end{cases}
\]
which by construction solves the Moser equation $$\partial_u g_{s, u} +dg_{s, u}(X_{s,u}) = 0.$$

We have $\|X_{s, u}\| = \frac{|\partial_u g_{s,u}|}{\|\nabla_{w, v}g_{s, u}\| }\leq \frac{3|w|\|A\|_\infty}{\|\nabla_{w, v}g_{s, u}\|}$ and  $X_{s, u}=0$ where $\|v\| < R$.
Else,  for $\|v\|\geq R$, $|w|\leq \sqrt R$ we have
$$
\begin{aligned}
\|X_{s,u}\|
&\leq
\frac{|\partial_u g_{s,u}|}
{\|\nabla_{v,w}g_{s,u}\|}
\\ &\leq \frac{|\partial_u g_{s,u}|}{\|\partial_{v}g_{s,u}\|}
\\ &\leq \frac{|\partial_u g_{s,u}|}{ \|\partial_vQ(v)\| -\|\partial_v\rho_s\|_\infty -3\varepsilon |w|\, \|A\|_\infty\ \partial_v\chi_R\|_\infty }
\\ &\leq
\frac{3\varepsilon\sqrt R\,\|A\|_\infty}{\lambda R
-\|\partial_v\rho_s\|_\infty
-3\varepsilon\sqrt R\,
\|A\|_\infty\|\partial_v\chi_R\|_\infty}.
\end{aligned} $$
     
    and for $\|v\|\geq R$, $|w|\geq \sqrt{R}$, 
$$\begin{aligned}
\|X_{s,u}\|
&\leq
\frac{3\varepsilon |w|\,\|A\|_\infty}
{\|\nabla_{w, v}g_{s,u}\|}
\\
&\leq
\frac{3\varepsilon |w|\,\|A\|_\infty}
{|\partial_wg_{s,u}|}
\\
&\leq
\frac{3\varepsilon |w|\,\|A\|_\infty}
{3|w|^2
-3\varepsilon
-3\varepsilon\|A\|_\infty
-\|\partial_w\rho_s\|_\infty}
\\
&\leq
\frac{3\varepsilon\|A\|_\infty}
{3\sqrt R
-\frac1{\sqrt R}(
3\varepsilon+3\varepsilon\|A\|_\infty
+\|\partial_w\rho_s\|_\infty
)}.
\end{aligned}$$

    For  $R$ sufficiently large $X$ is uniformly bounded. Since $X_{s, u}$ is vertical, its time $1$ flow is a fibered diffeomorphism $\Phi_s$ satisfying $h_s\circ \Phi_s = g_{s, 0}$. 

\end{proof}

\begin{deff}
    If $h(x; v) : M\times\R^n\map \R$ is a generating family, define its \textbf{difference function} by  \[\delta h: (x, v, v') \in M\times \R^n\times \R^n\mapsto h(x; v)- h(x; v') \in \R.\]
    The critical points of  $\delta h$ of critical value $l$ are in bijection with the Reeb chords of length $l$ of $\iota_h(\Sigma_h).$
    If $\bar\iota_h(x, v) = \bar\iota_h(x, v')$ and $\textbf{im}(d\bar\iota_{h, (x, v)}) \pitchfork \textbf{im}(d\bar\iota_{h, (x, v')})$ then $\delta h$ is of Morse type at $(x, v, v')$.
    \end{deff}

\begin{prop}
  For $C$ sufficiently large, the fiber of $\delta h_0$ over $C$ is diffeomorphic to a vector bundle over $M.$  
\end{prop}
\begin{proof}
 We have that $$\delta h_0(x, w, v, w', v') =  \\w^3-w'^3 +3\ep(w-w')+3w\eta(x, v) -3w'\eta(x, v')+Q(v)-Q(v').$$  
    Fix $x\in M$ and $v, v' \in \R^k$. If $v\not \in \operatorname{supp}(\eta)$, we can write $$\delta h_0(x, w, v, w', v') = w^3+3\ep w +R(x, v, v', w')$$ for some function $R$. Since $w^3+3\ep w$ is strictly increasing, for all choices of $w'$ there is a unique  $w$ such that $\delta h_0(x, w,v, w', v')=C$, depending smoothly on the other variables.
    The same argument holds if $v'\not \in \operatorname{supp}(\eta)$.
    Else, if $(v, v')\in \operatorname{supp}(\eta)^2$ we have that $b := Q(v)-Q(v')$, $a := 3\ep+3\eta(x, v)$ and $a' := -3\ep-3\eta(x, v')$ are bounded and we are looking for the solutions to  $$w^3+aw-w'^3+a'w'+b-C = 0.$$
    Fix $d = w+w' \in \R$. The equation becomes $$w^3+aw-d^3+3d^2w-3dw^2+w^3+a'd-a'w +b-C =0$$ with discriminant $$\Delta_w
= -\bigl(3d^2+2(a-a')\bigr)^3
  -27\bigl((a+a')d+2(b-C)\bigr)^2. $$
If $3d^2 > 2|a-a'|$  then $\Delta_w < 0$.
Else we know that $ |d| < \sqrt{\frac 23 |a-a'|}$ is also bounded, so $\Delta_w\xrightarrow[]{C\map\infty}-\infty$ uniformly, and we can choose $C$ large enough so that $\Delta_w < 0$ uniformly in $a, a', b, d$. For each $d \in \R$ the equation thus has a unique solution $w$ that depends smoothly on the other variables.  

Let $p : M\times \R\times \R^k\times\R\times\R^k\map M\times \R^k\times\R^k$ be the obvious projection. For $(x, v, v')\in M\times \R^k\times \R^k$ and $U$ a small neighbourhood of this point, we showed that  $\delta h_0\inv(C)\cap p\inv(U)$ is a graph above  $U\times \R$. Therefore $\delta h_0\inv(C)$ is diffeomorphic to a vector bundle over  $M$.
\end{proof}

We denote by $\overline{M} = M\times\R\times\R^k\times\R\times\R^k$ the domain of the functions $\delta h_s$. We know that $n_s$, d efined earlier, is a regular value of $\delta h_s$ and we set $$r^- = \{(x, s)\in \overline M \times [0, 1] \,|\, \delta h_s(x) = n_s\}.$$
Fix $n^+\in \R$ a real number such that for $C\geq n^+$,  $\delta h_0\inv(C) $ is a vector bundle over $M$ and $C$ is a regular value of all the functions $\delta h_s$. Set $$r^+ = \{(x, s)\in \overline M \times [0, 1] \,|\, \delta h_s(x) = n^+\}.$$

\begin{prop}
    There exists a family of smooth vector fields $X_s$ on $\overline{M}$ such that $X_s\delta h_s \geq 1$ outside of a compact set and on $r^\pm.$
\end{prop}
\begin{proof}
   We are looking for a pseudogradient for $\delta h_s(x, w, v, w', v') = w^3+3\ep w-w'^3-3\ep w'+Q(v)-Q(v')+\rho_s(x, v, w) -\rho_s(x, v', w') + w\eta_s(x, v)-w'\eta_s(x, v')$. 
    Set $X _s= \partial_w-\dron{w'}+\nabla Q(v)-\nabla Q(v')$.
    Since $Q$ is non degenerate, there exists a strictly positive constant $C_1$ such that $C_1\|v\|\leq \|\nabla Q(v)\| $ for all $v\in \R^k.$ Moreover, $|X_s\bigl(\rho_s(x, v, w)-\rho_s(x, v', w')\bigr)|$ is of compact support and thus bounded by a constant $C_2 \in \R$, and one can find a constant $C_3\in \R$ such that $$|X_s  (w\eta(x, v)-w'\eta(x, v'))|\leq 2\|\eta\|_\infty+ (|w|+|w'|)C_3.$$
We have
\begin{align*}
         X_s\bigl( w^3+3\ep w-w'^3-3\ep w'+Q(v)-Q(v')\bigr)  &= 3w^2 +3\ep+3w'^2+3\ep+\|\nabla Q(v)\|^2+\|\nabla Q(v')\|^2 \\
         &\geq 6\ep +3w^2+3w'^2+C_1^2\|v\|^2 +C_1^2\|v'\|^2\\
         &\geq 6\ep +\min(C_1^2, 3)\|(w, v, w', v')\|^2 .
    \end{align*}
Consequently, for $x \not\in \operatorname{supp}(\rho_s)\cup \operatorname{supp}(\eta_s)$  :

$$X_s\bigl(\delta h_s\bigr) \geq 6\ep.$$

Setting $C_4 = \min(C_1^2, 3)$, 
    \begin{align*}
         X_s\bigl(\delta h_s\bigr) &\geq  6\ep +C_4\|(w, v, w', v')\|^2-|X_s\bigl(\rho_s(x, v, w)-\rho_s(x, v', w')\bigr)|-|X_s  (w\eta(x, v)-w'\eta(x, v'))| \\
          &\geq  6\ep +C_4\|(w, v, w', v')\|^2-C_2-2\|\eta\|_\infty -C_3|w| -C_3|w'| \\
         &\geq  6\ep +C_4\|(w, v, w', v')\|^2-C_2-2\|\eta\|_\infty -\frac{C_4}2w^2-\frac{C_3^2}{C_4^2}-\frac{C_4}2w'^2-\frac{C_3^2}{C_4^2} \\
         &\geq  6\ep +\frac{C_4}2\|(w, v, w', v')\|^2-C_2-2\|\eta\|_\infty -\frac{2C_3^2}{C_4^2}
    \end{align*}
  so that $X_s(\delta h_s) \geq 6\ep$  for $\|(w, v, w', v')\|$ sufficiently large. Therefore $X_s(\delta h_s)  \geq 6\ep$ outside of a compact set.

Replacing $X_s$ with $aX_s + a\nabla\delta h_s$ with $a$ sufficiently large,  we can suppose $\delta h_s(X_s) \geq 1 $ on $r^\pm$ and outside of a compact set.
\end{proof}

\begin{figure}[htbp]
\centering

\begin{minipage}[t]{0.48\textwidth}
\centering
\resizebox{\linewidth}{!}{%
\begin{tikzpicture}[line cap=round, line join=round]

  \filldraw[fill=gray!5, draw=gray!35]
    (-4.3,-2.6) rectangle (4.5,2.8);
  \node[anchor=north east] at (-2.4,2.75)
    {$\overline M\times\{s\}$};

  \path[fill=blue!10]
    (-2,1.6) -- (2,1.6)
    .. controls (2.7,0.6) and (2.5,-0.7) .. (2.2,-1.6)
    -- (-2.2,-1.6)
    .. controls (-2.5,-0.7) and (-2.7,0.6) .. (-2,1.6)
    -- cycle;

  \draw[gray!65] (-3.9,1.6) -- (3.9,1.6);
  \draw[gray!65] (-3.9,-1.6) -- (3.9,-1.6);
  \node[above] at (3.55,1.6) {$r_s^+$};
  \node[below] at (3.55,-1.6) {$r_s^-$};

  \path[fill=orange!35, draw=orange!75!black, thick]
    (-1.1,0.25)
    .. controls (-1.15,0.75) and (-0.45,0.8) .. (0,0.55)
    .. controls (0.5,0.85) and (1.15,0.55) .. (1.1,0.1)
    .. controls (1.3,-0.45) and (0.65,-0.85) .. (0.15,-0.65)
    .. controls (-0.5,-0.95) and (-1.3,-0.5) .. cycle;
  \node[orange!60!black] at (0,-0.1)
    {$\mathcal K'_s$};

  \draw[orange!75!black, densely dashed, ->]
    (-1.1,0.25)
    .. controls (-1.15,0.8) and (-1,1.2) .. (-0.9,1.6);
  \draw[orange!75!black, densely dashed, ->]
    (1.1,0.1)
    .. controls (1.2,0.7) and (1.1,1.2) .. (0.95,1.6);

  \draw[blue!65!black, very thick]
    (-2,1.6) -- (2,1.6);
  \draw[blue!65!black, very thick]
    (-2.2,-1.6) -- (2.2,-1.6);

  \draw[orange!85!black, line width=3pt]
    (-0.9,1.6) -- (0.95,1.6);

  \node[blue!65!black, above] at (-1.65,1.65)
    {$\mathcal K_s^+$};
  \node[orange!60!black, above] at (0,1.65)
    {$\mathcal K_s^{\prime+}$};
  \node[blue!65!black, below] at (0,-1.65)
    {$\mathcal K_s^-=\partial^-\mathcal K_s$};

  \draw[blue!65!black, very thick, ->]
    (-2.2,-1.6)
    .. controls (-2.5,-0.7) and (-2.7,0.6) .. (-2,1.6);
  \draw[blue!65!black, very thick, ->]
    (2.2,-1.6)
    .. controls (2.5,-0.7) and (2.7,0.6) .. (2,1.6);

  \node[blue!65!black, anchor=east] at (-2.85,0)
    {$\partial^{//}\mathcal K_s$};
  \draw[blue!65!black] (-2.8,0) -- (-2.425,0);

  \node[blue!65!black] at (1.7,-0.8)
    {$\mathcal K_s$};

  \node[blue!65!black, anchor=west] at (2.85,0.35)
    {$X_s$};
  \draw[blue!65!black] (2.8,0.35) -- (2.43,0.35);

\end{tikzpicture}%
}

\vspace{2mm}
\small Slice of $\mc K$ in \(\overline{M}\times \{s\}. \)
\end{minipage}
\hfill
\begin{minipage}[t]{0.48\textwidth}
\centering
\resizebox{\linewidth}{!}{%
\begin{tikzpicture}[line cap=round, line join=round]

  \filldraw[fill=gray!5, draw=gray!35]
    (-5,-2.2) rectangle (4.2,2.2);
  \node at (3.8,1.85) {$r^+$};

  \fill[blue!10] (-3,-1.35) rectangle (3,1.35);

  \path[fill=orange!35, draw=orange!75!black, thick]
    (-3,0.35)
    .. controls (-2.4,0.35) and (-2.1,0.65) .. (-1.5,0.65)
    .. controls (-0.6,0.65) and (0.3,0.90) .. (1,0.90)
    .. controls (1.8,0.90) and (2.4,0.45) .. (3,0.45)
    -- (3,-0.40)
    .. controls (2.3,-0.40) and (2,-0.65) .. (1.3,-0.65)
    .. controls (0.4,-0.65) and (-0.3,-0.90) .. (-1,-0.90)
    .. controls (-1.9,-0.90) and (-2.3,-0.35) .. (-3,-0.35)
    -- cycle;

  \draw[<-, violet!80!black, densely dashed]
    (-3,0.90) -- (1,0.90);
  \draw[<-, violet!80!black, densely dashed]
    (-3,-0.90) -- (-1,-0.90);

  \draw[blue!65!black, very thick]
    (-3,-1.35) rectangle (3,1.35);

  \draw[violet!80!black, line width=5pt]
    (-3,-0.90) -- (-3,0.90);

  \draw[orange!85!black, line width=2pt]
    (-3,-0.35) -- (-3,0.35);

  \node[blue!65!black, above] at (-3,1.35)
    {$\mathcal K_0^+$};
  \node[blue!65!black, above] at (3,1.35)
    {$\mathcal K_1^+$};

  \node[blue!65!black] at (1.9,1.12)
    {$\mathcal K^+$};
  \node[orange!60!black] at (0,0)
    {$\mathcal K^{\prime+}$};

  \node[violet!80!black, anchor=east] at (-3.5,0.75)
    {$\mathcal K_0^{\prime+}$};
  \draw[violet!80!black]
    (-3.45,0.75) -- (-3,0.65);

  \node[blue!65!black, above] at (0,1.65)
    {$\partial \mathcal K_s^+$};
  \draw[blue!65!black] (0,1.65) -- (0,1.35);

  \draw[->] (-3,-1.8) -- (3.4,-1.8) node[right] {$s$};
  \draw (-3,-1.75) -- (-3,-1.85) node[below] {$0$};
  \draw (3,-1.75) -- (3,-1.85) node[below] {$1$};

\end{tikzpicture}%
}

\vspace{2mm}
\small Slice of $\mc K^+$ in $r^+.$
\end{minipage}
\label{fig:tranche-K}
\end{figure}
\begin{prop}
    There is a cobordism with boundaries $\mc K \subset \overline M\times [0, 1]$ such that
    \begin{enumerate}
        \item  The interior of $\mc K$ contains all the points $p=(x, s)$ such that $d(\delta h_s)_p(X_s)<1$ and $\delta h_s(x) \in [n^-_s, n^+]$ simultaneously.
        \item $\partial \mc K$ can be decomposed as $\partial^{//} \mc K \cup \partial^+\mc K\cup \partial^-\mc K $  
        \item $X_s$ is tangent to $\partial^{//}\mc K $.
        \item  $\partial^+ \mc K \subset r^+$ and $\partial^-\mc K \subset r^-$. 
        \item $\mc K_s := \mc K\cap \overline{M}\times \{s\}$ is an isotopy in $s$.
        \item $\partial^+ \mc K_0 \map  \delta h_0\inv(\{C\})$ is a homotopy equivalence.
    \end{enumerate}
\end{prop}
\FloatBarrier
\begin{proof}
Choose a compact set $ \mc K' \subset \overline{M} \times [0,1]$ satisfying $(1)$. By shrinking  $\mc K'$, we can suppose $\delta h_s(x) \in ]n^-_s, n^+[$  for$(x, s)\in \mc K'$.

Outside of $\mc K'$ define a vector field $Z = \partial_s- \frac{\partial_s\delta h_s-\dot n_s}{\|\nabla \delta h_s\|^2}\nabla\delta h_s$ and extend it arbitrarily to all of $\overline{M}\times [0, 1]$, with the condition $ds(Z) = 1$.  By construction the flow of $Z$ preserves $r^+$ and $r^-$.

Let  $\mc K'^+ = \varphi_X^{[0, C]}(\mc K') \cap r^+$ be the (compact) set of points of $r^+$ pushed into $\mc K'$ by the flow of $-X.$ 
Set $\mc K'^+_0= \{\fee_Z^{-s}(x, s)\, ;\, (x, s)\in \mc K'^+\} \subset \overline{M}\times \SingZ$, and choose a compact submanifold with boundary $\mc K_0^+$ of $r^+_0=\delta h_0\inv(n^+)\times \SingZ\subset \overline{M}\times [0, 1]$ whose interior contains $\mc K'^+_0$. Arrange that $\mc K_0^+\map r_0^+$ be a homotopy equivalence which will imply $(6)$. Set $\mc K_s^+ = \fee_Z^s(\mc K_0^+)$, which by construction satisfies $\mc K^+\supset \mc K'^+.$ Finally define $$\mc K= \{(x, s) \in \overline{M} \times [0, 1]\,|\, \delta h_s(x) \in [n_s^-, n^+] \text{ et } \exists u, \fee_X^u(x)\in \mc K^+\cup \mc K'\},$$
$$ \partial ^{//}\mc K :=  \{(x, s) \in \overline{M} \times [0, 1]\,|\, \delta h_s(x) \in [n_s^-, n^+] \text{ et } \exists u, \fee_X^u(x)\in \mc \partial K^+\}$$ and $$\partial^\pm\mc K = \mc K\cap r^\pm.$$

The trajectories of $-X_s$ that reach $\partial \mc K^+$ cannot intersect $\mc K'$, since by definition the trajectories reaching $\mc K'$  start in $\mc K'^+$ and $\mc K'^+\cap \mc  \partial\mc K^+ = \emptyset$. The union of these trajectories $\partial^{//}\mc K$ is then a smooth face of $\mc K$ diffeomorphic to $\partial\mc K^+ \times [0, 1]$.

We construct the isotopy showing $(5)$. For $x\in\partial\mc K_0^+$, let $T(x,s)$ denote the time taken by the trajectory of $-X_s$ starting from $\fee_Z^s(x) \in r_s^+$ to reach $r_s^-$. We have the following parametrisation of $\partial \mc K^{//}$, fibered above $[0,1]$ :

\[
\begin{aligned}
p:\partial\mc K_0^+\times[0,1]\times[0,1]
&\longrightarrow \partial^{//}\mc K,\\
(x,u,s)
&\longmapsto
\bigl(\fee_{X_s}^{-uT(x,s)}\circ\fee_Z^s(x),s\bigr).
\end{aligned}
\]
By construction of $p$, the vector field $V=p_*(\partial_s)$ is tangent to $\partial \mc K^{//}$, to  $\partial^{//}\mc K \cap \partial^+\mc K$ and to $\partial^{//}\mc K \cap \partial^-\mc K$, and satisfies $ds(V)=1$. We extend $V$ to the faces $\partial^\pm\mc K$ keeping these properties (for example by setting $V' = \tilde V-(\tilde V\delta h_s)\cdot\frac{X_s}{X_s\delta h_s}$ where $\tilde V$ is an arbitrary extension of $V$ to $\partial^\pm \mc K$), and then extend it to all of $\mc K$ . We obtain this way a vector field $V'$ tangent to the faces $\partial^{//}\mc K$ and $\partial^\pm\mc K$. Since moreover $ds(V') = 1$, its time $s$ flow is an isotopy
\[
(\mc K_0,\partial^-\mc K_0)
   \xrightarrow{\ \sim\ }
(\mc K_s,\partial^-\mc K_s),
\]
which shows $(5)$.

\end{proof}

\begin{deff}
    Let $(K, K^+, K^-)$ be a  cobordism. Define
\[
\mu(K, K^+, K^-)
=
\min_{\substack{f}}
\left\{
\#\operatorname{Crit}(f)
\;\middle|\;
f : K \to \mathbb{R}
\text{ of Morse type}
\right\}.
\]
   If $M$ is a smooth manifold, define \[
\mu_{\mathrm{st}}(M)
=
\min_{\substack{p>0\\ f}}
\left\{
\#\operatorname{Crit}(f)
\;\middle|\;
f : M \times \mathbb{R}^p \to \mathbb{R}
\text{ of Morse type and quadratic at infinity}
\right\}.
\]
   
\end{deff}

Let $\mc D(M)$  denote the set of free $\Lambda$-complexes $D_\bullet$ such that
\[
D_i = 0 \qquad \text{for } i<0 \text{ and } i>n,
\]
and belong to the simple equivalence class of  $C_\bullet(h)$, where $h: M\map \R$ is a function of Morse type.

\cite{DAMIAN_2002} shows that
\begin{prop}
\[
\mu_{\mathrm{st}}(M)
=
\min\left\{
\operatorname{rk}(D_\bullet)
\;\middle|\;
D_\bullet\in \mc D(M)
\right\}.
\]
\end{prop}

\begin{lemm}
    Let $(K, K^+, K^-)$ be a cobordism with boundary faces and $f : K\map \R$ a Morse-Bott function on  $K$ with only one critical set $L$, and a trivial positive normal bundle.
    Suppose furthermore that the inclusion $L \map K$ induces a left invertible map $\pi_1(L)\map \pi_1(K).$
    Then $$\mu(K) \geq \mu_{st}(L).$$
\end{lemm}
\begin{proof}
   Set $\Lambda_K = \ZZ[\pi_1(K)]$ and $\Lambda_L = \ZZ[\pi_1(L)]$.
    Let $h : L \map \R$ be a Morse function on $L$, Morse-Smale relative to a metric $\|\cdot\|_L$ and $l : K \map \R$ a Morse function on $K$ with $$\operatorname{rk}(\mc C(l, \Lambda_K)) =\mu(K).$$
    First perturb $f$ to a Morse function $\tilde f : K\map \R$ satisfying some properties. 
    Choose Morse-Bott coordinates around $L$ such that $ f(x, v) = Q(v)$ for $Q : \R^{2k} \map \R$ a quadratic form. By taking a cutoff function $ a  : \R \map \R$ with $\|a'\|_\infty$ sufficiently small, we can set $$\tilde f(x, v) = Q(v) +a(\|v\|^2) h(x),$$ which is of Morse type. By construction $\mc C(\tilde f, \Lambda_K)_{\bullet+k}  = \mc C(h, \Lambda_L)_{\bullet}\otimes \Lambda_K$, and so $$\mc C(\tilde f, \Lambda_K)_{\bullet+k} \otimes \Lambda_L  = \mc C(h, \Lambda_L)_{\bullet}\otimes \Lambda_K\otimes \Lambda_L\approx \mc C(h, \Lambda_L)_{\bullet}.$$
     
Moreover  $\mc C(\tilde f, \Lambda_K) \simeq_s\mc C(l, \Lambda_K)$ entails :  $$C(h, \Lambda_L)_{\bullet-k}\simeq_s\mc C(\tilde f, \Lambda_K) \otimes \Lambda_L\simeq_s\mc C(l, \Lambda_K)\otimes \Lambda_L.$$
Since $\mc C(l, \Lambda_K)$ is a free complex, $\operatorname{rk}(\mc C(l, \Lambda_K)\otimes \Lambda_L) = \operatorname{rk}(\mc C(l, \Lambda_K)) = \mu(K)$. By definition, $\mu_{\mathrm{st}}(L) = 
\min \left\{
\operatorname{rk}(D_\bullet)
\;\middle|\;
D_\bullet \simeq_s \mc C(h, \Lambda_L)
\right\}$ which proves the desired inequality.

\end{proof}
\begin{findepreuve}
    By hypothesis $\mc K_0^+ \approx W$, and topologically, $\mc K_0$ is the attachment of an index $N$ Morse-Bott handle to $\mc K_0^+$, so taking $N\geq 3$, $\pi_1(\mc K_0) \approx \pi_1(W)$. If $L_0\map W $ induces a left invertible map $\pi_1(L_0)\map \pi_1(W)$, it is also the case for$L_0 \approx \textbf{Crit}(\delta g_0) \map \mc K_0$, so we can apply  lemma $(21)$.

In the end, using the Morse inequality : 
$$\begin{aligned}
    \mu_{st}(L) &\leq \mu (\mc K_0, \mc K_0^+, \mc K_0^-) \\ &= \mu (\mc K_1, \mc K_1^+, \mc K_1^-) \\ &\leq \#\textbf{Crit}(\delta h_{1|{\mc K_1}}) \\ &=  \#(\pi(L_0)\cap \pi(L_1)).
    \end{aligned}
    $$

The function $\delta h_0$ is Morse-Bott with a unique critical set $\pi(L_0)\subset M$ of critical value in $]n_0, n^+[$. Let $D$ denote its descending disc : $$D = \{\textbf{x} \in \overline{M}\,|\,\lim_t\phi^t_{X_0}(\textbf{x})\in L \text{ et } \delta h_0(\textbf{x})\geq n_0\}.$$
The pair $(\mc K_0, \mc K_0^-)$ is equivalent to $(D, \partial D)$. The manifold $D$ can be identified with the unit disc bundle of the negative definite vector bundle $E$ of $\delta h_0$ at $L_0$. Since $E$ is a stabilisation of the trivial bundle $N_{L_0}(M)$, the bundle $E$ is itself trivial. In particular $E$ is an orientable vector bundle. Let $T(E)$ denote the Thom space associated to $E$, then 
\begin{align*}
H^*(\mc K_0, \mc K_0^-)&\approx H^*(D, \partial D) \\& \approx  H^*(T(E), *)\\&\approx H^*(T(E), T(E)\setminus L_0) \\&\approx H^{*-n}(L_0)
\end{align*}
  using the Thom isomorphism. By the universal coefficient theorem it follows that : $$\displaystyle\sum_i b_i(L_0) = \displaystyle\sum_i b_i(\mc K_0, \mc K_0^-) = \displaystyle\sum_ib_i(\mc K_1, \mc K_1^-).$$ When $\pi(L_0)\pitchfork \pi(L_1)$, the function $\delta h_1$ is of Morse type on $(\mc K_1,\mc  K_1^-)$ and the Morse inequality implies : $$\# \textbf{Crit}(\delta h_{1|\mc K_1})\geq \sum_i b_i(\mc K_1, \mc K_1^-).$$

In the end we conclude : 
\[
\begin{aligned}
\sum_i b_i(L) 
&\leq  \#\mathrm{Crit}(\delta h_ {1|\mc K_1}) \\
&= \#\pi(L_0)\cap \pi(L_1).
\end{aligned}
\]

\end{findepreuve}

\printbibliography
 \end{document}